\documentclass{article}
\usepackage{amsmath,amssymb}

\newtheorem{theorem}{Theorem}
\newtheorem{lemma}{Lemma}
\newtheorem{conjecture}{Conjecture}
\newtheorem{example}{Example}

\begin{document}

\title{Proof of a conjectured q,t-Schr\"{o}der identity}

\author{William J. Keith \\
Drexel University, 3141 Chestnut, Philadelphia, PA 19104, USA \\
wjk26@drexel.edu}

\maketitle

\begin{abstract}  A conjecture of Chunwei Song on a limiting case of the $q,t$-Schr\"{o}der theorem is proved combinatorially.  The proof matches pairs of tableaux to Catalan words in a manner that preserves differences in the maj statistic.
\end{abstract}

\section{Introduction}

In \cite{Song}, Chunwei Song examined the combinatorics of a limiting case of the $q,t$-Schr\"{o}der Theorem, proved by Haglund in \cite{Haglund}. To state the main conjecture which concerns this paper, use the standard notation $${(a;q)}_n = (1-a)(1-aq)\dots (1-aq^{n-1}), \left[ n \atop k \right] = \frac{{(q;q)}_n}{{(q;q)}_k{(q;q)}_{n-k}}$$ and let $S_{n,i}$ be the set of permutations on $n$ letters with longest increasing subsequence of length $i$.  Define \begin{multline*} F_{n,k} := \sum_{i=1}^k \sum_{{a_1 + \dots + a_k = n} \atop {a_i > 0}} q^{\sum_{i=1}^k \left( a_i \atop 2 \right)} t^{\sum_{i=1}^{k-1} (k-i) a_i} \frac{1}{{(t^k;q)}_{a_1} {(q,q)}_{a_k}} \\ \times \prod_{i=1}^{k-1} \left[ {{a_i + a_{i+1} - 1} \atop a_i} \right] \frac{1}{({t^{k-i};q)}_{a_i + a_{i+1}}} \times {(q;q)}_n {(t;t)}_n \end{multline*} and, with maj the major index of a permutation, $$H_{n,k} (q,t) := \sum_{i=1}^k \sum_{\sigma \in S_{n,i}} q^{\text{maj}(\sigma)} t^{\left( {n \atop 2} \right) - \text{maj}\left(\sigma^{-1}\right)} \, \text{.}$$

Then Song conjectured, for the $q=t$ case, 

\begin{conjecture} $F_{n,k}(q, q) = H_{n,k} (q, q) \, \text{.} $\end{conjecture}

For $k=1$, both sides are $q^{\left( {n \atop 2} \right)}$.  For $k>1$, the conjecture remains.  Song rewrote the $k=2$ case as (Conjecture 3.7 in \cite{Song}):

\begin{conjecture} $$\sum_{w \in CW_n} q^{\text{maj}(w) - n \, \text{des}(w)} = \sum_{\sigma \in S_n(123)} q^{\text{maj}(\sigma) - \text{maj}(\sigma^{-1})} \, \text{.}$$
\end{conjecture}

Here $S_n (123)$ is the set of 123-avoiding permutations (those with longest increasing subsequence of length no more than 2), and $CW_n$ is the set of \emph{Catalan words}: words $x_1 \dots x_{2n}$ with $x_i \in \{0,1\}$ such that there are exactly $n$ 1s and no initial segment has more 1s than 0s.  The \emph{descent set} of such a word is $Des(w) = \{i \vert x_i = 1, x_{i+1} = 0 \}$.  The maj of such a word is $\sum_{x_i \in Des(w)} x_i$ and $\text{des}(w) = \vert Des(w) \vert$. Song proved a subset of this conjecture, proving equality for the set of $w$ and $\sigma$ representing involutions, where $w$ is an involution if $w = w^{-1}$ under the inversion operation that reverses $w$ and exchanges 0 and 1, an operation that corresponds to reversing the associated Dyck path. This article proves the full $k=2$ case by a bijection that may be generalizable to any $k$.

\section{The bijection}

Since $\text{maj}(w^{-1}) = 2n \, \text{des}(w) - \text{maj}(w)$ (from the proof of Theorem 3.10 in \cite{Song}), we can rewrite Conjecture 2 more symmetrically as 

\begin{theorem}\label{MainTheorem} $$\sum_{w \in CW_n} q^{\text{maj}(w) -  \text{maj}(w^{-1})} = \sum_{\sigma \in S_n(123)} q^{2(\text{maj}(\sigma) - \text{maj}(\sigma^{-1}))} \, \text{.}$$
\end{theorem}

A map that takes Catalan words of length $2n$ to permutations $\sigma \in S_n(123)$ by a mapping $\phi(w)=\sigma$ such that $\text{maj}(w) -  \text{maj}(w^{-1}) = 2(\text{maj}(\sigma) - \text{maj}(\sigma^{-1}))$ would now prove the theorem.  This we construct.

Define a \emph{Catalan half-word} of length $n$ to be a 01 word $x_1 \dots x_n$ with no initial segment having more 1s than 0s.  The number of 1s may vary, though the initial-segment condition means there are never more than $\lfloor \frac{n}{2} \rfloor$.  The set of ordered pairs $(w_1, w_2)$ of Catalan half-words of length $n$ in which $w_1$ and $w_2$ have the same number of 1s is in bijection with $CW_n$ by $(w_1, w_2) \rightarrow w_1 w_2^{-1}$, inversion being the same as before (reversal and exchanging 0 and 1) and juxtaposition denoting simple concatenation.  Clearly $(w_2, w_1)$ corresponds to the inverse of $(w_1, w_2)$.

A \emph{tableaux}, or standard Young tableaux, of $n$ boxes is an arrangement of the numbers 1 through $n$ on points of the lattice in the first quadrant such that rows increase when read left to right and columns increase when read from the bottom up, and there are no gaps between occupied lattice points in rows or columns.  The set of occupied lattice points is the shape of such a tableaux.  The descent set $Des(\tau)$ of a tableaux $\tau$ is the set of $i$ with $i+1$ in a higher row than $i$ in $\tau$, and maj($\tau$) is the sum of such $i$.  Pairs of tableaux $(\tau_1,\tau_2)$ of $n$ boxes wherein both have the same shape are in bijection with permutations $ \sigma \in S_n$ through the famous Robinson-Schensted correspondence, which has the properties that pairs $(\tau_2, \tau_1)$ correspond to $\sigma^{-1}$ (so that diagonal pairs $(\tau_1,\tau_1)$ correspond to involutions), tableaux with exactly $i$ columns correspond to permutations in $S_{n,i}$, and maj($\tau_2$) = maj($\sigma$).

\begin{lemma}The set $T_D$ of tableaux $\tau$ of $n$ boxes, with no more than 2 columns and $j \geq 0$ entries in the second column, that have $D = Des(\tau) = \{i_1, \dots , i_k\}$, is equal in cardinality to the set $CW_D$ of Catalan half-words $w$ of length $n$ that have $j$ 1s and patterns $x_ix_{i+1}=01$ where $n-i \not\in D$.
\end{lemma}

\begin{example}  The two tableaux with descent set $D=\{1,3,4,6 \}$, $n=7$, $j=3$ are 

\begin{center}\begin{tabular}{c}
\framebox[0.75cm]{7\phantom{4}}\phantom{\framebox[0.75cm]{74}} \\
\framebox[0.75cm]{5\phantom{4}}\framebox[0.75cm]{6\phantom{4}} \\
\framebox[0.75cm]{2\phantom{4}}\framebox[0.75cm]{4\phantom{4}} \\
\framebox[0.75cm]{1\phantom{4}}\framebox[0.75cm]{3\phantom{4}} \\
\end{tabular} = $\tau_1$, 
\begin{tabular}{c}
\framebox[0.75cm]{5\phantom{4}}\phantom{\framebox[0.75cm]{54}} \\
\framebox[0.75cm]{4\phantom{4}}\framebox[0.75cm]{7\phantom{4}} \\
\framebox[0.75cm]{2\phantom{4}}\framebox[0.75cm]{6\phantom{4}} \\
\framebox[0.75cm]{1\phantom{4}}\framebox[0.75cm]{3\phantom{4}} \\
\end{tabular} = $\tau_2$.\end{center}

The permutation $\sigma$ associated to the pair $(\tau_1,\tau_2)$ by the Robinson-Schensted correspondence is $(7562143)$, with maj 14, and its inverse is $(5476231)$, also with maj 14, so their difference is 0.  The two words in $CW_{\{2,5 \}}$ for $n=7$, $j=3$ are $w_1 = 0010011$ and $w_2 = 0011010$.  The Catalan word $w$ associated to the ordered pair $(w_1, w_2)$ is $w=00100111010011$.  Then $w^{-1} = 00110100011011$, and $\text{maj}(w) - \text{maj}(w^{-1}) = 3+8+10-4-6-11 = 0$.
\end{example}

\noindent \textbf{Proof of Lemma 1}

Select $n$, $j$, and $D=\{i_1,\dots,i_k\}$.  These choices fix the following parts of the tableaux and Catalan half-words under consideration:

\[ \overbrace{0 \cdots 00 \underbrace{01}_{n-i_k,n-i_k+1} 11 \cdots100 \cdots 0   \underbrace{01}_{n-i_{k-1},n-i_{k-1}+1}  \cdots \underbrace{01}_{n-i_1,n-i_1+1} 11 \cdots100 \cdots 0 }^n \]

\begin{center}\begin{tabular}{c}
\framebox[1.2cm]{$i_k$\phantom{$\vdots_1$}}\phantom{\framebox[1.2cm]{7$\vdots_1$}} \\
\framebox[1.2cm]{\vdots\phantom{$\vdots_1$}}\phantom{\framebox[1.2cm]{7$\vdots_1$}} \\
\framebox[1.2cm]{$i_2$ \phantom{$\vdots_1$}}\framebox[1.2cm]{\vdots \phantom{$\vdots_1$}} \\
\framebox[1.2cm]{\vdots \phantom{$\vdots_1$}}\framebox[1.2cm]{$i_k+1$\phantom{$\vdots_1$}} \\
\framebox[1.2cm]{$i_1$ \phantom{$\vdots_1$}}\framebox[1.2cm]{\vdots\phantom{$\vdots_1$}} \\
\framebox[1.2cm]{\vdots \phantom{$\vdots_1$}}\framebox[1.2cm]{$i_2+1$\phantom{$\vdots_1$}} \\
\framebox[1.2cm]{2\phantom{$\vdots_1$}}\framebox[1.2cm]{\vdots \phantom{$\vdots_1$}} \\
\framebox[1.2cm]{1\phantom{$\vdots_1$}}\framebox[1.2cm]{$i_1+1$\phantom{$\vdots_1$}} \\
\end{tabular}\end{center}

In the word, the segments between 01 entries must be filled with nonincreasing sequences: $00\cdots 0$, $11\cdots10\cdots0$, or $11\cdots1$.  The choice of how to fill each segment is thus completely determined by the number of 1s we choose to enter.  This choice is bounded by the initial-segment condition, the number of spaces in each segment, and the total number of 1s available ($j-k$).

In the tableaux, the sections between $i_a+1$ and $i_{a+1}$ must be filled without producing any more descents.  The only way to do this is by listing the entries from $i_a+2$ in the right-hand column up to some $i_a+c$, and then listing the remaining entries from $i_a+c+1$ to $i_{a+1}-1$ in the left-hand column.  The choice of how to fill each section is thus completely determined by the choice of $c$.  This choice is bounded by the condition that the two parts of the fill cannot overlap, the length of the interval between $i_a$ and $i_{a+1}$, and the total number of spaces available in the right-hand column ($j-k$).

The two sets of conditions are exactly analogous.  We can fill the tableaux section between $i_1+1$ and $i_2$ by choosing a segment of length up to $i_1-1$ entries in the right-hand column, unless bounded by the spaces available.  We can fill the word segment between the last descent and the end of the word with a number of 1s up to the space available, $i_1 - 1$ spaces, unless bounded by the number of 1s available.  Fill each segment up the tableaux and back along the word, making equal legal choices at each step.$\Box$

This correspondence defines our map.  Let $\phi: S_n (123) \rightarrow CW_n$ by sending $\sigma \rightarrow (\tau_1, \tau_2) \rightarrow (w_1, w_2)$, where the first map is the Robinson-Schensted correspondence and the second sends the tableaux $\tau_i$ to the word $w_i$ such that $w^{-1}$ has the same descent set as $\tau_i$ and the number of 1s in the segments between 01 entries at positions $n-i_{a-1}$ and $n-i_a$ are given by the number of elements in the right-hand column of $\tau_1$ above each $i_a+1$.

\begin{example} The tableaux $\tau_1, \tau_2$ and words $w_1, w_2$ from Example 1 are from corresponding $T_D$ and $CW_D$.  The fixed and selectable sections in the words of $CW_D$ look like  $[0][01][?][01][?]$.  Since in $\tau_1$ the entry 4 appears above 3 in the right-hand column, the map assigns $\phi(\tau_1) = 0010011 = w_1$, filling the segment after position $n-i_1$ with a 1.  Since in $\tau_2$ the entry 4 appears in the left-hand column, and the 7 appears in the right-hand column after $i_2+1=6$, the map fills the segment right of $n-i_2=n-5$ with the 1, giving $\phi(\tau_2)= 0011010=w_2$.
\end{example}

All tableaux with the same descent set have the same maj.  (Other tableaux may have the same maj statistic; there are as many such sets as there are partitions of maj($w$) into distinct nonconsecutive parts less than $n$.)  As can be seen above, not all of the words in $CW_D$ have the same maj; instead, it will be shown in the next section that each makes the same contribution to maj($w$) - maj($w^{-1}$) under the correspondence of Catalan half-words with Catalan words described above.  To complete the proof of the theorem it remains only to show that maj($\tau_i$) and maj($w_i$) are related as claimed.

\section{Proofs of required properties}

\begin{lemma} Let $w_1$ and $w_2$ be two words in some $CW_D$, thus having the same length $n$, the same number of 1s $j$, and the same locations of 01 patterns starting at $\{n-i_1,\dots,n-i_k\}$.  Let $w_3$ be any other Catalan half-word with $j$ 1s.  Say $w_b = w_1 w_3^{-1}$ and $w_c = w_2 w_3^{-1}$.  Then maj($w_b$) - maj($w_b^{-1}$) = maj($w_c$) - maj($w_c^{-1}$).
\end{lemma}

\noindent \textbf{Proof} If $w_1$ and $w_2$ are different, there will be an earliest digit in which they differ, say $x_i$, with $w_1$ having an entry of 0 and $w_2$ having an entry of 1.  Since by hypothesis they possess the same number of 1s, there will be a latest digit, say $x_j$, in which $w_1$ has an entry of 1 and $w_2$ has an entry of 0.  Suppose we alter $w_1$ by exchanging its 0 and 1 entry in these two places.  Call $w_4$ the resulting word.  (The initial-segment condition may be violated in $w_4$, but the maj statistic is defined for all 01-words.)  What is the relation between maj($w_1 w_3^{-1}$) - maj(${(w_1 w_3^{-1})}^{-1})$ and maj($w_4 w_3^{-1}$) - maj(${(w_4 w_3^{-1})}^{-1})$?

\begin{align*} w_b = w_1 w_3^{-1} &= 0\cdots x_{i-1}0x_{i+1}\cdots x_{j-1}1x_{j+1}\cdots w_3^{-1} \\
w_b^{-1} = w_3 w_1^{-1} &= w_3 \cdots (1-x_{j+1}) 0 (1-x_{j-1}) \cdots (1-x_{i+1}) 1 (1- x_{i-1}) \cdots 1 \\
w_4 w_3^{-1} &= 0\cdots x_{i-1}1x_{i+1}\cdots x_{j-1}0x_{j+1}\cdots w_3^{-1} \\
w_3 w_4^{-1} &= w_3 \cdots (1-x_{j+1}) 1 (1-x_{j-1}) \cdots (1-x_{i+1}) 0 (1- x_{i-1}) \cdots 1 
\end{align*}

Note that $x_{j+1}$ may be the first digit of $w_3^{-1}$.  Additionally, it does not matter to maj($w$) - maj($w^{-1}$) whether the transition between $w_1$ and $w_3^{-1}$ is a descent, or not; if it is, then $w_1$ ends in a 1 and the first digit of $w_3^{-1}$ is a 0, and so the last digit of $w_3$ is 1 and the first digit of $w_1$ is 0 as well.  Likewise, if the transition is not a descent in $w$, it is not in $w^{-1}$ either.

Because $w_1$ and $w_2$ are both in the same $CW_D$, we cannot have 01 patterns in different places in the two words.  This proscribes the case where $x_{i-1}=0$, since then $x_{i-1} x_i$ would be 00 in $w_1$, but would be a 01 sequence in $w_2$.  Similarly, we cannot have $x_{i+1} = 1$ (in particular, we cannot have $x_j = x_{i+1}$), $x_{j-1} = 0$, or $x_{j+1} = 1$ with the exception of the possibility that $x_{j+1}$ is the first digit of $w_3^{-1}$.  If it is not, then our situation is

\begin{align*} w_b = w_1 w_3^{-1} &= 0 \cdots 100 \cdots 110 \cdots w_3^{-1} \\
w_b^{-1} = w_3 w_1^{-1} &= w_3 \cdots 100 \cdots 110 \cdots 1 \\
w_4 w_3^{-1} &= 0 \cdots 110 \cdots 100 \cdots w_3^{-1} \\
w_3 w_4^{-1} &= w_3 \cdots 110 \cdots 100 \cdots 1
\end{align*}

\noindent in the case when $x_{j+1}$ is not the first digit of $w_3^{-1}$.  We have descents in $w_b$ at positions $x_{i-1}$ and $x_j$, and in $w_b^{-1}$ at positions $x_{2n-j}$ and $x_{2n-i+1}$, contributing to maj($w_b$) - maj($w_b^{-1}$) a total of $i-1+j-(2n-j+2n-i+1) = 2i+2j-4n-2$.  In $w_4 w_3^{-1}$, the first descent is moved forward a place, and the second is moved back 1, and likewise for the descents in the inverse.  Thus the total does not change.

If $x_{j+1}$ is the first digit of $w_3^{-1}$, then our situation is

\begin{align*} w_b = w_1 w_3^{-1} &= 0 \cdots 100 \cdots 11 \, w_3^{-1} \\
w_b^{-1} = w_3 w_1^{-1} &= w_3 \, 00 \cdots 110 \cdots 1 \\
w_4 w_3^{-1} &= 0 \cdots 110 \cdots 10 \,  w_3^{-1} \\
w_3 w_4^{-1} &= w_3 \, 10 \cdots 100 \cdots 1
\end{align*}

If $x_{j+1}=1$, then the last digit of $w_3$ is 0.  Then $w_b$ has a descent at $x_{i-1}$ and $w_b^{-1}$ a descent at $x_{2n-i+1}$, contributing a total of $2i-2n-2$.  Meanwhile $w_4 w_3^{-1}$ has descents at $x_i$ and $x_{n-1}$, and its inverse has descents at $x_{n+1}$ and $x_{2n-i}$, contributing a total of $i+n-1-(n+1+2n-i)=2i-2n-2$.  If $x_{j+1} = 0$, $w_b$ possesses another descent at $x_n$ and $w_b^{-1}$ does as well, cancelling out, whereas $w_4 w_3^{-1}$ possesses neither.

Since maj($w_b$) - maj($w_b^{-1}$) does not alter under this operation, and by proceeding down the two words performing such exchanges we will eventually cause them to match, any two words in the same $CW_D$ make the same contribution to maj($w$) - maj($w^{-1})$. $\Box$

\begin{lemma} If $w_1$ is in $CW_D$ for $D=\{i_1, \dots, i_k \}$, in which words have length $n$ and exactly $j$ 1s, then the contribution of $w_1$ to maj($w_1 w_2^{-1}$) - maj($w_2 w_1^{-1}$) is $2j - 2(maj(\tau_1))$ for $\tau_1$ any tableaux with descent set $D$.
\end{lemma}

\noindent \textbf{Proof} Since all elements of any $CW_D$ contribute the same amount to maj($w$)-maj($w^{-1}$), we can select one of the simpler ones and calculate what that contribution is.  Let us examine the word in which 1s are entered as far to the right as possible:

\begin{align*} {{[0 \cdots 0 ]} \atop {l_1}} {{01} \atop {}} & {{[0 \cdots 0 ]} \atop {l_2}} {{01} \atop {}} {\cdots \atop {}}  {{01} \atop {}} {{[1 \cdots 1 0 \cdots 0 ]} \atop {l_g}} {{01} \atop {}} {{[1 \cdots 1 ]} \atop {l_{g+1}}} {{01} \atop {}} {\cdots \atop {}} {{01} \atop {}} {{[1 \cdots 1 ]} \atop {l_{k+1}}} 
\end{align*}

It is notationally convenient to use $l_i$ to denote the lengths of the segments between 01 patterns rather than differences between the $i_k$ counting backward from $n$.  Note that any of these lengths may be 0, in which case $l_g$ may be nonunique and any nearby empty segment may be assigned the label.  If there is no segment in which 1s and 0s are mixed, the calculation that follows is similar.

What is the contribution this word makes to maj($w_1 w_2^{-1}$) - maj($w_2 w_1^{-1}$)?

In $w_1 w_2^{-1}$, the first descent occurs at place $l_1 + 2$.  The next occurs at place $l_1 + l_2 + 4$, and so forth up to $l_1 + \dots + l_{g-2} + 2(g-2)$.  The descent following the segment labeled $l_{g-1}$ is further ahead of this pattern by a distance equal to the number of 1s remaining after the $k$ 1s in the 01 patterns are accounted for and the first segments of length $l_{g+1}$ through $l_{k+1}$ are filled: thus it occurs at $l_1 + \dots + l_{g-1} + 2(g-1) + (j-k - (l_{g+1} + \dots + l_{k+1})$.  The descent following the segment labeled $l_g$ comes at the end of the following full segment of 1s, so it happens at place $l_1 + \dots + l_g + 2g + l_{g+1}$.  This pattern holds up to the descent following the segment labeled $l_{k-1}$, which occurs after the full segment of 1s labeled $l_k$, thus at place $l_1 + \dots + l_{k-1} + 2(k-1) + l_k$.  The segments labeled $l_k$ and $l_{k+1}$ do not have descents following them; as in the previous discussion, whether or not a descent occurs in the transition from the end of $w_1$ to the beginning of $w_2^{-1}$ does not matter, since if one does it will be cancelled upon subtraction.

The sum of the places of descents in the $w_1$ half of $w$ is thus

\[ l_1 (k-1) + l_2 (k-2) + \dots + l_{k-1} (1) + (k-1)(k) + (j-k) - l_{k+1} \, \text{.} \]

In $w^{-1}$, the inverted form of $w_1$ appears after $w_2$, thus: 

\begin{align*} {{(w_2)} \atop {}}  {{[0 \cdots 0 ]} \atop {l_{k+1}}} {{01} \atop {}}  {{[0 \cdots 0 ]} \atop {l_k}} {{01} \atop {}} {\cdots \atop {}}  {{01} \atop {}} {{[1 \cdots 1 0 \cdots 0 ]} \atop {l_g}} {{01} \atop {}} {{[1 \cdots 1 ]} \atop {l_{g-1}}} {{01} \atop {}} {\cdots \atop {}} {{01} \atop {}} {{[1 \cdots 1 ]} \atop {l_{1}}} 
\end{align*}

\noindent and we note that $w_2$ is of length $n = l_1 + \dots + l_{k+1} + 2k$.

Here the first descent in this half of the word comes at place $n + l_{k+1} +2$.  The next occurs at $n+l_{k+1}+l_k+4$, and so forth until we reach $l_g$.  The number of 1s in this segment is the number of 0s in the segment's original state.  Thus the relevant descent occurs at $n + l_{k+1} + \dots + l_{g+1} + 2(k-g+1) + l_g - (j-k-(l_{g+1} + \dots + l_{k+1}))$.  The next descent occurs at place $n+l_{k+1} + \dots + l_g + 2(k-g+2) + l_{g-1}$, and so forth until the last descent occurs at place $n + l_{k+1} + \dots + l_3 + 2(k-1) + l_2$.  Totaling, we find a contribution to maj($w^{-1}$) of

\begin{multline*} (n)(k-1) + (k-1)(k) + l_{k+1} (k) + \dots + l_2 (1) - (j-k) \\
 = (l_1 + \dots l_{k+1} + 2k)(k-1) + (k-1)(k) + l_{k+1} (k) + \dots + l_2 (1) - (j-k) \, \text{.}
\end{multline*}

Subtract the latter contribution from the former: the total contribution of $w_1$ to maj($w$) - maj($w^{-1}$) is

\begin{multline*} 2(j-k) - 2k(k-1) + l_1(0) + l_2 (-2) + l_3 (-4) + \cdots + l_{k+1} (-2k) \\
= 2j - 2k^2 + l_1(0) + l_2 (-2) + l_3 (-4) + \cdots + l_{k+1} (-2k) \, \text{.}
\end{multline*}

For the associated set of tableaux, it is straightforward to total up in terms of the $l_i$ the maj statistic of a tableaux with descent set identified by this $D$, confirming it to be $l_{k+1} k + l_{k} (k-1) + \cdots + l_2 (1) + k^2$.  Thus the contribution of $w_1$ to maj($w$) - maj($w^{-1}$) is $2j - 2(\text{maj}(\tau_1))$, as required.

$\Box$

Now we have that $w_1$ contributes $2j - 2(\text{maj}(\tau_1))$ to maj($w$) - maj($w^{-1}$), where maj($\tau_1$) = maj($\sigma^{-1}$).  Likewise $w_2$, in the opposite position, contributes $-(2j - 2(\text{maj}(\tau_2))$, so that the total of the two contributions is $2 (\text{maj}(\tau_2) -  \text{maj}(\tau_1))$: exactly twice the difference between the maj statistics of each tableaux, which are the maj statistics of the corresponding permutations $\sigma$ and $\sigma^{-1}$ respectively.  Summing over all Catalan half-words and corresponding tableaux, we have proved the theorem. $\Box$

\section{Generalization}

The next step would seem to be applying a similar map to 012 words and permutations in $S_n(1234)$ via tableaux of up to 3 columns.  We would need to find the appropriate rewrite of Conjecture 1, determine the necessary rule for the generalized Catalan word, and work out the correct map.  The larger goal, naturally, would be to prove Conjecture 1 fully by performing a similar procedure for the general case of $01\dots k$ words and tableaux of up to $k+1$ columns.

\end{document}